\documentclass[a4paper,12pt]{article}

\usepackage[T2A]{fontenc}
\usepackage[cp1251]{inputenc}
\usepackage{amsfonts,amsmath,amsthm,amscd,amssymb,latexsym,amsbsy,pb-diagram,cite,caption2}
\usepackage[dvips]{graphicx}
\usepackage[mathscr]{eucal}

\theoremstyle{plain}
\newtheorem{d_theorem}{Theorem}[section]

\newtheorem{d_proposition}{Proposition}[section]
\newtheorem{d_corollary}{Corollary}[section]
\newtheorem{d_definition}{Definition}[section]

\theoremstyle{definition}
\newtheorem{d_remark}{Remark}[section]
\newtheorem{d_example}{Example}[section]

\numberwithin{equation}{section}

\begin{document}
\begin{center}
{\Large\bf Tangent spaces to metric spaces and to their subspaces}
\end{center}

\bigskip
\begin{center}
{\bf O. Dovgoshey}
\end{center}

\begin{abstract}
We investigate a tangent space at a point of a general metric
space and metric space valued derivatives.  The condi\-tions under
which two different subspace of a metric space have isometric
tangent spaces in a common point of these subspaces are completely
deter\-minated.
\end{abstract}



\section{Introduction. Tangent metric spaces}

The recent achievements in the metric space theory are closely
related to some generalizations of the differentiation. The
concept of the upper gradient  \cite{He,HeKo,Sh}, Cheeger's notion
of differentiability for Rademacher's theorem in certain metric
measure spaces \cite{Ch}, the metric derivative in the studies of
metric space valued functions of bounded variation \cite{Am,AmTi}
and the Lipshitz type approach in \cite{Ha} are interesting and
important examples of such generalizations. A very interesting
technical tool to develop a theory of a differentiation in metric
separable spaces is the fact that every separable metric space
admits an isometric embedding into the dual space of a separable
Banach space. It provides a linear structure, and so a
differentiation, for a separable metric space, see for example a
rather complete theory of rectifiable sets and currents on metric
spaces in \cite{AmKi1,AmKi2}.

These generalizations of the differentiability usually lead to
nontrivial  results only for the assumption that metric spaces
have ``sufficiently many'' rectifiable curves. In almost all
mentioned approaches we see that  theories of  differentiations in
metric spaces involve an induced linear structure that is able to
use the classical differentiations in the linear normed spaces.

A new, intrinic, notion of differentiabililty for the mappings
between the general metric spaces was produced by O.~Dovgoshey and
O.~Martio in \cite{DM}. A basic technical tool in \cite{DM} is a
tangent space to an arbitrary metric space $X$ at a point $a\in X$
that was defined as a factor space of a family of sequences of
points $x_{n}\in X$ \ which converge to $a.$ This approach makes
possible to define a metric space valued derivative of  functions
$f:X\rightarrow Y,$ $X$ and $Y$ are metric spaces, as a mapping
between tangent spaces to $X$ at the point $a$ and, respectively,
to $Y$ at the point $f(a)$. The analysis of general properties of
tangent spaces and of metric space valued derivatives is the main
purpose of the present paper.

Let $(X,d)$ be a metric space and let $a$ be point of $X$. Fix
some sequence $\tilde r$ of positive real numbers $r_n$ which tend
to zero. In what follows this sequence $\tilde r$ be called a
\textit{normalizing sequence}. Let us denote by  ${\tilde X}$ the
set of all sequences of points from $X$.

\begin{d_definition}
\label{1:d1.1} Two sequences $\tilde{x},\tilde{y}\in \tilde X$,
$\tilde x=\{ x_n\}_{n\in \mathbb{N}}$ and $\tilde y=\{ y_n\}_{n\in
\mathbb{N}}$,  are mutually stable (with respect to a normalizing
sequence $\tilde r=\{ r_n\}_{n\in \mathbb{N}}$) if there is a
finite limit
\begin{equation}  \label{1:eq1.1}
\lim_{n\to\infty} \frac{d(x_n, y_n)}{r_n} := \tilde d_{\tilde
r}(\tilde x, \tilde y)=\tilde d(\tilde x, \tilde y).
\end{equation}
\end{d_definition}

We shall say that a family $\tilde F\subseteq \tilde X$ is
\textit{self-stable} (w.r.t. a normalizing sequence $\tilde r$) if
every two $\tilde x,\tilde y \in \tilde F$ are mutually stable. A
family $\tilde F \subseteq\tilde X$ is \emph{maximal self-stable}
if $\tilde{F}$ is self-stable and for an arbitrary $\tilde z\in
\tilde X$ either $\tilde z\in \tilde F$ or there is $\tilde x\in
\tilde F$ such that $\tilde x$ and $\tilde z$ are not mutually
stable.

A standard application of Zorn's Lemma leads to the following

\begin{d_proposition}
\label{1:p1.2}
Let $(X, d)$ be a metric space and let $a\in X$.
Then for every normalizing sequence $\tilde r=\{ r_n\}_{n\in
\mathbb{N}}$ there exists a maximal self-stable family $\tilde
X_a=\tilde X_{a, \tilde r}$ such that $\tilde a:=\{ a, a,
\dots\}\in \tilde X_a$.
\end{d_proposition}

Note that the condition $\tilde a \in \tilde X_a$ implies the
equality
\begin{equation*}
\lim_{n\to\infty} d(x_n, a)=0
\end{equation*}
for every $\tilde x=\{ x_n\}_{n\in \mathbb{N}}$ which belongs to
$\tilde X_a$.

Consider a function $\tilde d : \tilde X_a\times \tilde X_a \to
\mathbb{R}$ where $\tilde d(\tilde x, \tilde y)=\tilde d_{\tilde
r}(\tilde x, \tilde y)$ is defined by \eqref{1:eq1.1}. Obviously,
$\tilde d$ is symmetric and nonnegative. Moreover, the triangle
inequality for $d$ implies
\begin{equation*}
\tilde d (\tilde x, \tilde y) \leq \tilde d(\tilde x, \tilde
z)+\tilde d(\tilde z, \tilde y)
\end{equation*}
for all $\tilde x, \tilde y, \tilde z$ from $\tilde X_a$. Hence
$(\tilde X_a, \tilde d)$ is a pseudometric space.

\begin{d_definition}
\label{1:d1.3} The pretangent space to the space $X$ at the point
$a$ w.r.t. a normalizing sequence $\tilde r$ is the metric
identification of the pseudometric space $(\tilde X_{a, \tilde r},
\tilde d)$.
\end{d_definition}

Since the notion of pretangent space is basic for the present
paper, we remaind this metric identification construction.

Define a relation $\sim$ on $\tilde X_a$ by $\tilde x \sim \tilde
y$ if and only if $\tilde d(\tilde x, \tilde y)=0$. Then $\sim $
is an equivalence relation. Let us denote by
$\Omega_{a}:=\Omega_{a, \tilde r}=\Omega_{a,\tilde r}^X$ the set
of equivalence classes in $\tilde X_a$ under the equivalence
relation $\sim$. It follows from general properties of
pseudometric spaces, see, for example, \cite[Chapter 4,
Th.~15]{Kell}, that if $\rho$ is defined on $\Omega_a$ by
\begin{equation}  \label{1:eq1.2}
\rho(\alpha, \beta) :=\tilde d(\tilde x, \tilde y)
\end{equation}
for $\tilde x\in \alpha$ and $\tilde y \in \beta$, then $\rho$ is
the well-defined metric on $\Omega_a$. The metric identification
of $(\tilde X_a, \tilde d)$ is, by definition, the metric space
$(\Omega_a, \rho)$.

Remark that $\Omega_{a, \tilde r} \neq\emptyset$  because the
constant sequence $\tilde a$ belongs to $\tilde X_{a,\tilde r}$,
see Proposition \ref{1:p1.2}.

Let $\{n_k\}_{k\in \mathbb{N}}$ be an infinite, strictly
increasing sequence of natural numbers. Let us denote by $\tilde
r^{\prime}$ the subsequence $\{ r_{n_k}\}_{k\in\mathbb{N}}$ of the
normalizing sequence $\tilde r=\{ r_n\}_{n\in \mathbb{N}}$ and let
$\tilde x^{\prime}:=\{ x_{n_k}\}_{k\in \mathbb{N}}$ for every
$\tilde x=\{x_n\}_{n\in \mathbb{N}}\in\tilde{X}$.  It is clear
that if $\tilde x$ and $\tilde y $ are mutually stable w.r.t.
$\tilde r$, then $\tilde x^{\prime}$ and $\tilde y^{\prime}$ are
mutually stable w.r.t. $\tilde r^{\prime}$ and
\begin{equation}  \label{1:eq1.3}
\tilde d_{\tilde r} (\tilde x, \tilde y)=\tilde d_{\tilde
r^{\prime}}( \tilde x^{\prime},\tilde y^{\prime}).
\end{equation}
If $\tilde X_{a, \tilde r}$ is a maximal self-stable (w.r.t.
$\tilde r$) family, then, by Zorn's Lemma, there exists a maximal
self-stable (w.r.t. $\tilde r^{\prime}$) family $\tilde X_{a,
\tilde r^{\prime}}$ such that
\begin{equation*}
\{ \tilde x^{\prime}:\tilde x \in \tilde X_{a, \tilde
r}\}\subseteq \tilde X_{a, \tilde r^{\prime}}.
\end{equation*}
Denote by $\mathrm{in}_{\tilde r^{\prime}}$ the mapping from
$\tilde X_{a, \tilde r}$ to $\tilde X_{a, \tilde r^{\prime}}$ with
$\mathrm{in}_{\tilde r^{\prime}}(\tilde x)=\tilde x^{\prime}$ for
all $\tilde x\in \tilde X_{a, \tilde r}$. If follows from
\eqref{1:eq1.2} that after metric identifications
$\mathrm{in}_{\tilde r^{\prime}}$ pass to an isometric embedding
$\mathrm{em}^{\prime}$: $\Omega_{a, \tilde r}\to \Omega_{a, \tilde
r^{\prime}}$ under which the diagram
\begin{equation}  \label{1:eq1.4}
\begin{array}{ccc}
\tilde X_{a, \tilde r} & \xrightarrow{\ \ \mbox{in}_{\tilde r'}\ \
} &
\tilde X_{a, \tilde r^{\prime}} \\
\!\! \!\! \!\! \!\! \! p\Bigg\downarrow &  & \! \!\Bigg\downarrow
p^{\prime}
\\
\Omega_{a, \tilde r} & \xrightarrow{\ \ \mbox{em}'\ \ \ } &
\Omega_{a, \tilde r^{\prime}}
\end{array}
\end{equation}
is commutative. Here $p$, $p^{\prime}$ are metric identification
mappings, $p(\tilde x):=\{\tilde y\in \tilde X_{a, \tilde
r}:\tilde d_{\tilde r}(\tilde x,\tilde y)=0\}$ and
$p^{\prime}(\tilde x):=\{ \tilde y\in\tilde X_{a,\tilde
r^{\prime}}:\tilde d_{\tilde r^{\prime}}(\tilde x,\tilde y)=0\}$.

Let $X$ and $Y$ be two metric spaces. Recall that a map $f:X\to Y$
is called an \textit{isometry} if $f$ is distance-preserving and
onto.

\begin{d_definition}
\label{1:d1.4} A pretangent $\Omega_{a, \tilde r}$ is tangent if
$\mathnormal{em}^{\prime}$: $\Omega_{a, \tilde r}\to\Omega_{a,
\tilde r^{\prime}}$ is an isometry for every $\tilde r^{\prime}$.
\end{d_definition}

To verify the correctness of this definition, we must prove that
if $\tilde{X}_{a,\tilde{r}^{\prime }}^{(1)}$ and  $\tilde{X}_{a,
\tilde{r}'}$ are two distinct maximal self-stable families such
that the inclusions
\begin{equation}
\tilde{X}_{a,\tilde{r}^{\prime}}\supseteq \left\{\tilde{x}^{\prime
}:\tilde{x}\in \tilde{X}_{a,\tilde{r}}\right\}\subseteq
\tilde{X}_{a,\tilde{r}^{\prime}}^{(1)}\label{1:eq1.5}
\end{equation}%
hold and $em^{\prime }:\Omega
_{a,\tilde{r}}\,\,\longrightarrow\,\, \Omega
_{a,\tilde{r}^{\prime}}$ is an isometry, then $em_{1}^{\prime
}:\Omega _{a,\tilde{r}}\,\,\longrightarrow\,\, \Omega^{(1)}
_{a,\tilde{r}^{\prime}}$  is also an isometry, where $\Omega
_{a,\tilde{r}^{\prime }}^{(1)}$ is the metric identification of
$\tilde{X}_{a,\tilde{r}'}^{(1)}$. Indeed, it is clear that if
$\tilde{x}=\left\{
x_{n}\right\}_{n\in\mathbb{N}}\in\tilde{X}_{a,\tilde{r}},~~\tilde{y}=\left\{
y_{k}\right\}_{k\in\mathbb{N}}\in \tilde{X}_{a,\tilde{r}^{\prime
}}$ and
\[
\underset{k\rightarrow \infty }{\lim
}\frac{d(y_{k},x_{n_{k}})}{r_{n_{k}}}=0,
\]%
then there is  $\tilde{z}\in \tilde{X}_{a,\tilde{r}}$ with  $
\tilde{z}^{\prime}=\tilde{y}$. Consequently, since $em^{\prime }$
is an isometry and diagram \eqref{1:eq1.4} is commutative, the mapping $%
in_{r^{\prime }}:\tilde{X}_{a,\tilde{r}}\longrightarrow \tilde{X}%
_{a,\tilde{r}^{\prime }}$ is surjective, i.e.,%
\[
\tilde{X}_{a,\tilde{r}^{\prime }}=\left\{\tilde{x}^{\prime }:%
\tilde{x}\in \tilde{X}_{a,\tilde{r}}\right\} .
\]%
Hence, by \eqref{1:eq1.5}, we obtain the inclusion%
$
\tilde{X}_{a,\tilde{r}^{\prime }}^{(1)}\supseteq \tilde{X}_{a,
\tilde{r}^{\prime}}.
$
It implies the equality $\tilde{X}_{a,\tilde{r}^{\prime }}^{(1)}=%
\tilde{X}_{a,\tilde{r}^{\prime }}$ because $\tilde X_{a,\tilde
r'}$ is maximal self-stable. Hence $em_{1}^{\prime }=em^{\prime }$
and, so $em_{1}^{\prime }$ is an isometry.

These arguments give the following proposition.

\begin{d_proposition}
\label{1:p1.5} Let $X$ be a metric space with a marked point $a$,
$ \tilde{r}$ a normalizing sequence and $\tilde{X}_{a,\tilde{r}}$
a maximal self-stable family with correspondent pretangent space
$\Omega _{a,\tilde{r}}$. The following statements are equivalent.
\begin{itemize}
\item[$(i)$] $\Omega_{a,\tilde{r}}$ is tangent.

\item[$(ii)$] For every subsequence $\tilde{r}^{\prime }$ of the
sequence $ \tilde{r}$ the family $\left\{ \tilde{x}^{\prime
}:\tilde{x}\in \tilde{X}_{a,\tilde{r}}\right\}$ is maximal
self-stable w.r.t. $\tilde{r}^{\prime}$.

\item[$(iii)$] A function $em^{\prime
}:\Omega_{a,\tilde{r}}\longrightarrow \Omega _{a,\tilde{r}^{\prime
}}$ is surjective for every $\tilde{r} ^{\prime}$.

\item[$(iv)$] A function $in_r':\,\tilde{X}_{a,\tilde r}
\longrightarrow \tilde{X}_{a,\tilde{r}^{\prime}}$ is surjective
for every $\tilde{r}'$.
\end{itemize}
\end{d_proposition}

Now we introduce an equivalence relation for the classification of
normalizing sequences.

\begin{d_definition}
\label{1:d1.8} Let $X$ be a metric space with a marked point $a$.
Two normalizing sequences $\tilde{r}$ and $\tilde{t}$ are
equivalent at the point $a$,  $\tilde{r}$ $\approx $ $\tilde{t}$,
if the logical equivalence
$$
(\tilde{F} \text{ is self-stable w.r.t. }
\tilde{r})\Longleftrightarrow (\tilde{F} \text{ is self-stable
w.r.t. } \tilde{t})
$$
 is true for every $\tilde{F}\subseteq
\tilde{X}$ with $\tilde{a}\in \tilde{F}$.
\end{d_definition}

A normalizing sequence $\tilde{r}$ will be called
\textit{confluented in a point} $a$ if there exists an one-point
pretangent space $\Omega _{a,\tilde{r}}$  (it certainly implies
that all pretangent $\Omega _{a,\tilde{r}}$ are one-point).

\begin{d_theorem}
\label{1:p1.9} Let $(X,d)$ be a metric space with a marked point
$a$ and let  $\tilde{r}=\left\{ r_{n}\right\} _{n\in\mathbb{N}}$
and $\tilde{t}=\left\{ t_{n}\right\} _{n\in\mathbb{N}}$ be two
normalizing sequences which are equivalent at the point $a$. Then
at least one of the following statements holds.
\newline
$(i)$ There is a real number $c>0$ such that
\begin{equation}
\underset{n\rightarrow \infty }{\lim }\frac{r_{n}}{t_{n}}=c.
\label{1:eq1.7}
\end{equation}
\newline
$(ii)$ The sequences $\tilde{r}$ and $\tilde{t}$ are confluented
in the point $a$.
\end{d_theorem}
\begin{proof}
Suppose that both sequences $\tilde r$ and $\tilde t$ are not
confluented in  $a$. Then there are $\tilde{x}=\left\{
x_{n}\right\} _{n\in\mathbb{N}}$ and $\tilde{y}=\left\{
y_{n}\right\} _{n\in\mathbb{N}}$ from $\tilde{X}$ such that
\begin{equation}
\tilde{d}_{\tilde{r}}(\tilde{x},\tilde{a})=\underset{ n\rightarrow
\infty }{\lim }\frac{d(x_{n},a)}{r_{n}}>0 ~~ \text{and}
~~\tilde{d}_{
\tilde{t}}(\tilde{y},\tilde{a})=\underset{n\rightarrow \infty }{
\lim }\frac{d(y_{n},a)}{t_{n}}>0  \label{1:eq1.8}
\end{equation}%
where $\tilde{a}=(a,a,\dots)$. If $\tilde{d}_{\tilde{r}}(\tilde{
y},\tilde{a})>0$ or $~\tilde{d}_{\tilde{t}}(\tilde{x},
\tilde{a})>0$, then we obtain
\begin{equation*}
0<\frac{\tilde{d}_{\tilde{r}}(\tilde{y},\tilde{a})}{
\tilde{d}_{\tilde{t}}(\tilde{y},\tilde{a})}=\underset{n\rightarrow
\infty}{\lim }\frac{t_{n}}{r_{n}}<\infty
\end{equation*}%
or, respectively,
\begin{equation*}
0<\frac{\tilde{d}_{\tilde{t}}(\tilde{x},\tilde{a})}{
\tilde{d}_{\tilde{r}}(\tilde{x},\tilde{a})}=\underset{n\rightarrow
\infty}{\lim }\frac{r_{n}}{t_{n}}<\infty ,
\end{equation*}%
i.e., Statement $(i)$ holds. Now observe that the equalities
\begin{equation}
\tilde{d}_{\tilde{r}}(\tilde{y},\tilde{a})=\tilde{d}_{%
\tilde{t}}(\tilde{x},\tilde{a})=0  \label{1:eq1.9}
\end{equation}%
lead to a contradiction because \eqref{1:eq1.8} and \eqref{1:eq1.9} imply%
\begin{equation*}
0=\underset{n\rightarrow \infty }{\lim }\frac{r_{n}}{t_{n}}=\infty
.
\end{equation*}
Thus if Statement $(i)$ does not hold, then at least one of the
sequences $\tilde r$ and $\tilde t$ is confluented. We claim that
if  $\tilde r$ or $\tilde t$ is confluented, then both $\tilde r$
and $\tilde t$ are confluented. Indeed, if $\tilde r$ confluented
and we have a finite limit
\begin{equation}\label{1:eq1.9*}
\tilde d_{\tilde t}(\tilde y,\tilde
a)=\lim_{n\to\infty}\frac{d(y_n,a)}{t_n}\ne0
\end{equation}
for  $\tilde y=\{y_n\}_{n\in\mathbb N}\in\tilde X$. Then $\tilde
d_{\tilde r}(\tilde y,\tilde a)=0$ because $\tilde t\approx\tilde
r$ and $\tilde r$ is confluented in the point $a$. Write
\begin{equation}\label{1:eq1.10}
y_n^*:=\begin{cases} y_n&\text{if }n\text{ is odd}\\
a&\text{if $n$ is even}
\end{cases}
\end{equation}
for every $n\in\mathbb N$ and put $\tilde
y^*:=\{y_n^*\}_{n\in\mathbb N}$. Then we obtain $\tilde d_{\tilde
r}(\tilde y^*,\tilde a)=0$. Thus the family
$$
\tilde F:=\{\tilde y,\tilde y^*,\tilde a\}
$$
is self-stable w.r.t. $\tilde r$. Since $\tilde r\approx\tilde t$,
this family also is self-stable w.r.t $\tilde t$. Consequently
there is a finite limit
$$
\tilde d_{\tilde t}(\tilde a,\tilde
y^*)=\lim_{n\to\infty}\frac{d(y_n^*,a)}{t_n}.
$$
Hence, by \eqref{1:eq1.9*} and \eqref{1:eq1.10}, we obtain
$$
0\ne\lim_{n\to\infty}\frac{d(y_{2n+1}^*,a)}{t_{2n+1}}=\lim_{n\to\infty}\frac
{d(y_{2n}^*,a)}{t_{2n}}=0.
$$
This contradiction shows that $\tilde t$ is confluented if $\tilde
r$ is confluented.

Hence Statement $(i)$ holds if Statement $(ii)$ does not hold, and
the theorem follows.
\end{proof}

\begin{d_remark}\label{1:r1.8}It is clear that if there is $c>0$ such
that \eqref{1:eq1.7} holds, then normalizing sequences $\tilde{r}$
and $\tilde{t}$ are equivalent at every point $a\in X$.
\end{d_remark}
\begin{d_proposition}
\label{1:p1.10} Let $(X,d)$ be a metric space with a marked point
$a$. The following propositions are equivalent.
\begin{itemize}
\item[$(i)$] The point $a$ is an isolated point of the metric
space $X$.

\item[$(ii)$] Every two normalizing  sequences are equivalent at
the point $a$.

\item[$(iii)$] All normalizing sequences are confluented in $a$.
\end{itemize}
\end{d_proposition}

\begin{proof}
 The implication $(i)\Rightarrow(ii)$ and $(i)\Rightarrow(iii)$ are trivial. To
prove $(ii) \Rightarrow(i)$ suppose that the relation
\begin{equation}
\tilde{r}\approx\tilde{t}  \label{1:eq1.10*}
\end{equation}%
holds for every two normalizing $\tilde{r}$ and $\tilde{t}$ but
there is $\tilde{x}=\left\{ x_{n}\right\}_{n\in\mathbb{N}}\in
\tilde X$ such that $\underset{n\rightarrow \infty }{\lim
}d(x_{n},a)=0$ and $ d(x_{n},a)>0$ for all $n\in\mathbb{N}$. Let
$\tilde{x}^{\prime }=\left\{ x_{n_{k}}\right\}_{k\in\mathbb{N}}$
be an infinite subsequence of $\tilde{x}$ with
\begin{equation}
\underset{k\rightarrow \infty }{\lim
}\frac{d(x_{n_{k}},a)}{d(x_{k},a)}=0. \label{1:eq1.11}
\end{equation}%
Write%
\begin{equation*}
\tilde{r}:=\left\{ d(x_{k},a)\right\}
_{k\in\mathbb{N}},\qquad\tilde{t}:=\left\{ d(x_{n_{k}},a)\right\}
_{k\in\mathbb{N}}.
\end{equation*}%
Now, by the construction, both $\tilde{r}$ and $\tilde{t}$ \ are
not confluented and, moreover, \eqref{1:eq1.11} imply that
\eqref{1:eq1.7} does not hold for any $c>0$.  Hence, by Theorem
\ref{1:p1.9}, $\tilde{r}$ and $\tilde{t} $  are not equivalent at
the point $a$, contrary to \eqref{1:eq1.10*}. Thus the implication
$(ii)\Rightarrow (i)$ is true.

If $a$ is not an isolated point of $X$, then there is a sequence
$\tilde b=\{ b_n\}_{n\in N}\in\tilde X$ such that
$\lim\limits_{n\to\infty} d(a, b_n)=0$ and $d(a, b_n)\neq 0$ for
all $n\in \mathbb N$. Consider the normalizing sequence $\tilde
r=\{ r_n\}_{n\in N}$ with $r_n := d(a, b_n)$. It follows
immediately from \eqref{1:eq1.1} that $\tilde d_{\tilde r} (\tilde
a, \tilde b)=1$ where $\tilde a$ is the constant sequence $\{ a,
a, \dots\}$. The application of Zorn's Lemma shows that there is a
maximal self-stable family $\tilde X_{a, \tilde r}$ such that
$\tilde a, \tilde b \in \tilde X_{a, \tilde r}$. Then the metric
identification of the pseudometric space $(\tilde X_{a, \tilde r},
\tilde d)$ has at least two points. Consequently we also have
$(iii)\Rightarrow(i)$.
\end{proof}

\section{Metric space valued derivatives. \\Definition and general
properties}

Let $(X_{i},d_i),\ i=1,2$, be metric spaces with marked points
$a_{i}\in X_{i}$ and $\tilde{r}_{i}=\{r_n^{(i)}\}_{n\in\mathbb N}$
normalizing sequences and $\tilde{X}_{a_{i},\tilde{r}_{i}}^{i}$
maximal self-stable families with correspondent pretangent spaces
$\Omega_{a_{i},\tilde{r}_{i}}$. For  functions $f:X_{1}\rightarrow
X_{2}$ define the mappings $\tilde{f}:\tilde{X}_{1}\rightarrow
\tilde{X}_{2}$ as
\begin{equation}
\tilde{f}(\tilde{x})=\left\{ f(x_{i})\right\} _{i\in \mathbb{N}
}\text{ for  }\tilde{x}=\left\{ x_{i}\right\} _{i\in \mathbb{N}
}\in \tilde{X}_{1}.  \label{2:eq3.1}
\end{equation}

\begin{d_definition}
\label{2:d3.1} A function $f:X_{1}\rightarrow X_{2}$ is
differentiable w.r.t. the pair $\left(\tilde{X}_{a_{1},\tilde{r}_{1}}^{1},%
\tilde{X}_{a_{2},\tilde{r}_{2}}^{2}\right) $ if the following
conditions are satisfied:
\begin{itemize}
\item[$(i)$] $\tilde{f}(\tilde{x})\in \tilde{X}_{a_{2},\tilde{r}
_{2}}^{2}$ for every $\tilde{x}\in
\tilde{X}_{a_{1},\tilde{r}_{1}}^{1};$

\item[$(ii)$] $\left(  \tilde{d}_{\tilde{r}_{1}}(\tilde{x},
\tilde{y}) =0\right) \Longrightarrow \left(
\tilde{d}_{\tilde{r}_{2}}(\tilde f(\tilde{x}),\tilde
f(\tilde{y}))=0
\right) $ for all $\tilde{x},\tilde{y}\in \tilde{X}%
_{a_{1},\tilde{r}_{1}}^{1}$, where
$$
\tilde d_{\tilde r_1}(\tilde x,\tilde
y)=\lim_{n\to\infty}\frac{d_1(x_n,y_n)}{r_n^{(1)}},\qquad \tilde
d_{\tilde r_2}(\tilde f(\tilde x),\tilde f(\tilde
y))=\lim_{n\to\infty}\frac{d_2(f(x_n),f(y_n))}{r_n^{(2)}}.
$$
\end{itemize}
\end{d_definition}
\begin{d_remark}\label{r:2.1}
Note that condition $(i)$  of Definition \ref{2:d3.1} implies the
equality
$$
f(a_1)=a_2.
$$
\end{d_remark}
 Let $p_{i}:\tilde{X}_{a_{i},\tilde{r}_{i}}^{i}\rightarrow
\Omega _{a_{i},\tilde{r}_{i}},~i=1,2,$ be metric identification
mappings.

\begin{d_definition}
\label{2:d3.2} A function $D^{\ast }f:\Omega _{a_{1},\tilde{r}
_{1}}\rightarrow \Omega _{a_{2},\tilde{r}_{2}}$ is a metric space
valued derivative of  $f:X_{1}\rightarrow X_{2}$ at the point
$a_{1}\in X_{1}$ w.r.t. the pair $\left( \tilde{\Omega }_{a_{1},
\tilde{r}_{1}},\tilde{\Omega }_{a_{2},\tilde{r}_{2}}\right) $ (or,
in short, a derivative of $f$ ) if $f$ is differentiable w.r.t.
$\left( \tilde{X}_{a_{1},\tilde{r}
_{1}}^{1},\tilde{X}_{a_{2},\tilde{r}_{2}}^{2}\right) $ and the
following diagram
\begin{equation}
\begin{diagram}
\node{\tilde{X}_{a_{1},\tilde{r}_{1}}^{1}}
                        \arrow[2]{e,t}{\tilde{f}}
                        \arrow{s,l}{p_{1}}
\node[2]{\tilde{X}_{a_{2},\tilde{r}_{2}}^{2}}
                        \arrow{s,r}{ p_{2}} \\
\node{\tilde{\Omega }_{a_{1},\tilde{r}_{1}}}
\arrow[2]{e,t}{D^{\ast }f} \node[2]{\tilde{\Omega
}_{a_{2},\tilde{r}_{2}}}
\end{diagram}
\label{2:eq3.2}
\end{equation}%
is commutative.
\end{d_definition}

In this section we establish some common properties of the metric
space valued derivatives.

Let us show, first of all, that the metric space valued derivative
is unique if exists. Indeed, suppose that diagram \eqref{2:eq3.2}
is commutative with $D^{\ast }f=D_{1}^{\ast }f$ and with $D^{\ast
}f=D_{2}^{\ast }f$. Let $\beta \in \tilde{\Omega
}_{a_{1},\tilde{r}_{1}}^{1}.$ Since $p_{1}$ is a surjection, there
is $\tilde{x}_{1}\in \tilde{X}_{a_{1},\tilde{r}_{1}}^{1}$ such
that $ \beta =p_{1}(\tilde{x}_{1}). $ Definition \ref{2:d3.1}
implies that if $ \beta =p_{1}(\tilde{y}_{1}) $ for some other
$\tilde{y}_{1}\in \tilde{X}_{a_{1},\tilde{r} _{1}}^{1},$ then
\[
p_{2}(\tilde{f}(\tilde{x}_{1}))=p_{2}(\tilde{f}(\tilde{y}_{1})).
\]%
From the commutativity of \eqref{2:eq3.2} we obtain%
\[
D_{1}^{\ast }(\beta )=D_{1}^{\ast }(p_{1}(\tilde{x}_{1}))=p_{2}(
\tilde{f}(\tilde{x}_{1})) =D_{2}^{\ast
}(p_{1}(\tilde{x}_{1}))=D_{2}^{\ast }(\beta ),
\]%
i.e., $D_{1}^{\ast }=D_{2}^{\ast}$.

The following proposition shows that the Chain Rule remains valid
for the metric space valued derivatives.

\begin{d_proposition}
\label{2:p3.3} Let $X_{i}$ be metric spaces with marked points
$a_{i}\in X_{i}$ and $\tilde{r}_{i}$ normalizing sequences and
$\tilde{X}_{a_{i},\tilde{r}_{i}}^{i}$ maximal self-stable families
with
correspondent pretangent spaces $\Omega _{a_{i},\tilde{r}%
_{i}},~i=1,2,3.$

Let  $f:X_{1}\rightarrow X_{2}$ $\ $and  $g:X_{2}\rightarrow
X_{3}$ be differentiable functions, $f$  w.r.t. the pair $\left(
\tilde{X}_{a_{1},\tilde{r}_{1}}^{1},\tilde{X}_{a_{2},\tilde{r}
_{2}}^{2}\right) $ and $g$ w.r.t. $\left(\tilde{X}_{a_{2},
\tilde{r}_{2}}^{2},\tilde{X}_{a_{3},\tilde{r}_{3}}^{3}\right) .$

Then the superposition $\psi =g\circ f$ is differentiable w.r.t.
$\left(\tilde{X}_{a_{1},\tilde{r}_{1}}^{1},\tilde{X}_{a_{3},
\tilde{r}_{3}}^{3}\right)$ and
\begin{equation}
D^{\ast}(\psi )=(D^{\ast }g)\circ(D^{\ast}f).  \label{2:eq3.3}
\end{equation}
\end{d_proposition}

\begin{proof}
The differentiability of $\psi$ is an immediate consequence of the
differentiability  of  $f$ and $g,$ see Definition \ref{2:d3.1}.
To prove \eqref{2:eq3.3}  note that
\[
p_{2}\circ \tilde{f}=(D^{\ast }f)\circ p_{1}\text{ and }p_{3}\circ
\tilde{g}=(D^{\ast }g)\circ p_{2},
\]
see the following diagram
\begin{equation}\label{2:eq3.4}
\begin{diagram}
\node{\tilde X^1_{a_1,\tilde r_1}}  \arrow[2]{e,t}{\tilde \psi}
                                    \arrow{se,t}{\tilde f}
                                    \arrow[3]{s,l}{p_1}
\node[2]{\tilde X^3_{a_3,\tilde r_3}} \arrow[3]{s,r}{p_3} \\
\node[2]{\tilde X^2_{a_2,\tilde r_2}} \arrow{ne,t}{\tilde g}
                                      \arrow{s,r}{p_2}\\
\node[2]{\Omega_{a_2,\tilde r_2}}      \arrow{se,t}{D^*g} \\
\node{\Omega_{a_1,\tilde r_1}}     \arrow{ne,t}{D^*f}
                                 \arrow[2]{e,t}{D^*\psi}
\node[2]{\Omega_{a_3,\tilde r_3}}
\end{diagram}
\end{equation}
 Consequently we have
\begin{multline*}
p_{3}(\tilde{\psi })=p_{3}\circ (\tilde{g}\circ \tilde{f}
)=(p_{3}\circ \tilde{g})\circ \tilde{f}=((D^{\ast }g)\circ
p_{2})\circ\tilde{f}=
\\
=(D^{\ast }g)\circ (p_{2}\circ \tilde{f})=(D^{\ast }g)\circ
(D^{\ast }f)\circ p_{1},
\end{multline*}%
that is
\[
p_{3}\circ \tilde{\psi }=(D^{\ast }g\circ D^{\ast }f)p_{1}.
\]%
Hence the diagram
\begin{equation*}
\begin{diagram}
\node{\tilde{X}_{a_{1},\tilde{r}_{1}}^{1}}
\arrow[2]{e,t}{\tilde{\psi }} \arrow{s,l}{p_{1}}
\node[2]{\tilde{X}_{a_{2},\tilde{r}_{3}}^{3}}
\arrow{s,r}{p_{3}}\\
\node{{\Omega }_{a_{1},\tilde{r}_{1}}} \arrow[2]{e,t}{(D^{\ast
}g)\circ (D^{\ast}f)} \node[2]{{\Omega }_{a_{3},\tilde{r}_{3}}}
\end{diagram}
\end{equation*}
is commutative. The uniqueness of the derivative $D^{\ast }\psi $
and Definition \ref{2:d3.2} imply \eqref{2:eq3.3}.
\end{proof}
\begin{d_proposition}\label{d:p2.4}
Let $(X,d)$ and $(Y,\rho)$ be metric spaces, $a\in X$ and $b\in Y$
marked points in these spaces, and $f:X\to Y$  a function such
that $f(a)=b$. If for every maximal self-stable family $\tilde
X_{a,\tilde r}\subseteq\tilde X$ there is a maximal self-stable
family $\tilde Y_{b,\tilde t}\subseteq\tilde Y$ such that $f$ is
differentiable w.r.t. the pair $(\tilde X_{a,\tilde r},\tilde
Y_{b,\tilde t})$, then $f$ is continuous at the point $a$.
\end{d_proposition}
\begin{proof}
We may suppose that $a$ is not an isolated point of $X$. Let
$\tilde x=\{x_n\}_{n\in\mathbb N}\in\tilde X$ be a sequence such
that
$$\lim_{n\to\infty}d(x_n,a)=0\quad\text{and}\quad d(x_n,a)\ne0$$
for every $n\in\mathbb N$.
 Then, by Zorn's Lemma, there is a
maximal self-stable family
$$
\tilde X_{a,\tilde r}\supseteq\{\tilde a,\tilde x\}
$$
where $\tilde r=\{r_n\}_{n\in\mathbb N}$ is a normalizing sequence
with $r_n:=d(x_n,a)$ for  $n\in \mathbb N$. Hence there exists a
normalizing sequence $\tilde t=\{t_n\}_{n\in\mathbb N}$ for which
the  limit
$$
\lim_{n\to\infty}\frac{\rho(f(x_n),b)}{t_n}
$$
is finite. Consequently we have
$\lim_{n\to\infty}\rho(f(x_n),b)=0$ because
$$\lim_{n\to\infty}t_n=0.$$ Hence the function $f$ is continuous at
the point $a$.
\end{proof}

\section{Tangent spaces to subspaces of metric spaces}

Let $(X,d)$ be a metric space  with a marked point $a$, let $Y$
and $Z$ be subspaces of $X$ such that $a\in Y\cap Z$ and let
$\tilde r=\{r_n\}_{n\in\mathbb N}$ be a normalizing sequence.
\begin{d_definition}\label{6:d5.1}
The subspaces $Y$ and $Z$ are {\it tangent equivalent} at the
point $a$ w.r.t. the normalizing sequence $\tilde r$ if for every
$\tilde y_1=\{y_n^{(1)}\}_{n\in\mathbb N}\in\tilde Y$ and every
$\tilde z_1=\{z_n^{(1)}\}_{n\in\mathbb N}\in\tilde Z$ with finite
limits
$$
\tilde d_{\tilde r}(\tilde a, \tilde
y_1)=\lim_{n\to\infty}\frac{d(y_n^{(1)},a)}{r_n}\quad\text{and}\quad
\tilde d_{\tilde r}(\tilde a, \tilde
z_1)=\lim_{n\to\infty}\frac{d(z_n^{(1)},a)}{r_n}
$$
there exist $\tilde y_2=\{y_n^{(2)}\}_{n\in\mathbb N}\in\tilde Y$
and $\tilde z_2=\{z_n^{(2)}\}_{n\in\mathbb N}\in\tilde Z$ such
that
$$
\lim_{n\to\infty}\frac{d(y_n^{(1)},z_n^{(2)})}{r_n}=\lim_{n\to\infty}
\frac{d(y_n^{(2)},z_n^{(1)})}{r_n}=0.
$$
\end{d_definition}
We shall say that $Y$ and $Z$ are {\it strongly tangent
equivalent} at $a$ if $Y$ and $Z$ are tangent equivalent at $a$
for all normalizing sequences $\tilde r$.

Let  $\tilde F\subseteq\tilde X$. For a normalizing sequence
$\tilde r$ we define a family $[\tilde F]_Y=[\tilde F]_{Y,\tilde
r}$ by the rule
\begin{equation}\label{6:eq5.1}
(\tilde y\in[\tilde F]_Y)\Leftrightarrow((\tilde y\in\tilde
Y)\&(\exists\,\tilde x\in\tilde F:\tilde d_{\tilde r}(\tilde
x,\tilde y)=0)).
\end{equation}
Note that $[\tilde F]_Y$ can be empty for some nonvoid families
$\tilde F$ if the set $X\setminus Y$ is ``big enough''.
\begin{d_proposition}\label{6:p5.2}
Let $Y$ and $Z$ be subspaces of a metric space $X$ and let $\tilde
r$ be a normalizing sequence. Suppose that $Y$ and $Z$ are tangent
equivalent (w.r.t. $\tilde r$) at a point $a\in Y\cap Z$. Then
following statements hold for every maximal self-stable (in
$\tilde Z$) family $\tilde Z_{a,\tilde r}$.
\begin{itemize}
\item[$(i)$] The family $[\tilde Z_{a,\tilde r}]_Y$ is maximal
self-stable (in $\tilde Y$) and we have the equalities
\begin{equation}\label{6:eq5.2}
[[\tilde Z_{a,\tilde r}]_Y]_Z=\tilde Z_{a,\tilde r}=[\tilde
Z_{a,\tilde r}]_Z.
\end{equation}

\item[$(ii)$] If $\Omega^Z_{a,\tilde r}$ and $\Omega^Y_{a,\tilde
r}$ are metric identifications of $\tilde Z_{a,\tilde r}$ and,
respectively, of $\tilde Y_{a,\tilde r}:=[\tilde Z_{a,\tilde
r}]_Y$, then the mapping
\begin{equation}\label{6:eq5.3}
\Omega_{a,\tilde
r}^Z\ni\alpha\longmapsto[\alpha]_Y\in\Omega_{a,\tilde r}^Y
\end{equation}
is an isometry. Furthermore if $\Omega_{a,\tilde
r}^Z$ is tangent, then $\Omega^Y_{a,\tilde r}$ also is tangent.
\end{itemize}
\end{d_proposition}
\begin{proof}
$(i)$ Let $\tilde y_1,\tilde y_2\in\tilde Y_{a,\tilde r}:=[\tilde
Z_{a,\tilde r}]_Y$. Then, by \eqref{6:eq5.1}, there exist $\tilde
z_1,\tilde z_2\in \tilde Z_{a,\tilde r}$ such that
\begin{equation}\label{6:eq5.4}
\tilde d_{\tilde r}(\tilde y_1,\tilde z_1)=\tilde d(\tilde
y_2,\tilde z_2)=0.
\end{equation}
Since $\tilde z_1$ and $\tilde z_2$ are mutually stable, $\tilde
y_1$ and $\tilde y_2$ also are mutually stable. Consequently
$\tilde Y_{a,\tilde r}$ is self-stable. The similar arguments show
that $[\tilde Y_{a,\tilde r}]_Z$ is also self-stable. Moreover
since
$$
[[\tilde Z_{a,\tilde r}]_Y]_Z=[\tilde Y_{a,\tilde
r}]_Z\subseteq\tilde Z_{a,\tilde r},
$$
the maximality of $\tilde Z_{a,\tilde r}$ implies the first
equality in \eqref{6:eq5.2}. The second one also simply follows
from the maximality of $\tilde Z_{a,\tilde r}$. It still remains
to prove that $\tilde Y_{a,\tilde r}$ is a maximal self-stable
subset of $\tilde Y$. Let $\tilde Y^m_{a,\tilde r}$ be a maximal
self-stable family in $\tilde Y$ such that $\tilde Y^m_{a,\tilde
r}\supseteq\tilde Y_{a,\tilde r}$. Then $[\tilde Y^m_{a,\tilde
r}]_Z$ is  self-stable and $[\tilde Y^m_{a,\tilde
r}]_Z\supseteq\tilde Z_{a,\tilde r}$.  Since $\tilde Z_{a,\tilde
r}$ is maximal self-stable, the last inclusion implies the
equality $[\tilde Y^m_{a,\tilde r}]_Z=\tilde Z_{a,\tilde r}$.
Using this equality and \eqref{6:eq5.2} we obtain
$$
\tilde Y^m_{a,\tilde r}=[[\tilde Y^m_{a,\tilde r}]_Z]_Y=[\tilde
Z_{a,\tilde r}]_Y:=\tilde Y_{a,\tilde r},
$$
 i.e., $\tilde
Y_{a,\tilde r}$ is maximal self-stable.

$(ii)$ Let $\alpha\in\Omega^Z_{a,\tilde r}$ and let $\tilde
z\in\tilde Z_{a,\tilde r}$ such that $\tilde z\in\alpha$. It
follows from \eqref{6:eq5.1} that
\begin{equation}\label{6:eq5.5}
[\alpha]_Y=\{\tilde y\in\tilde Y:\tilde d_{\tilde r}(\tilde
y,\tilde z)=0\}.
\end{equation}
The last equality implies that function \eqref{6:eq5.3} is
distance-preserving. In addition, using \eqref{6:eq5.5} we see
that
$$
[[\alpha]_Y]_Z=\alpha\quad\text{and}\quad [[\beta]_Z]_Y=\beta
$$
for every $\alpha\in\Omega^Z_{a,\tilde z}$ and every
$\beta\in\Omega^Y_{z,\tilde r}$. Consequently function
\eqref{6:eq5.3} is bijective. To prove that $\Omega^Y_{a,\tilde
r}$ is tangent if $\Omega^Z_{a,\tilde r}$ is tangent we can use
Statement (ii) of Proposition~\ref{1:p1.5} and Statement (i) of
the present proposition.
\end{proof}
\begin{d_corollary}\label{6:c5.3}
Let $Y$ and $Z$ be  subspaces of a metric space $X$. Suppose that
$Y$ and $Z$ are tangent equivalent at a point $a\in Y\cap Z$
w.r.t. a normalizing sequence $\tilde r$ and that there exists a
unique maximal self-stable (in $\tilde Z$) family $\tilde
Z_{a,\tilde r}\ni\tilde a$. Then $\tilde Y_{a,\tilde r}:=[\tilde
Z_{a,\tilde r}]_Y$ is a unique maximal self-stable in $\tilde
Y_{a,\tilde r}$ which contains $\tilde a$.
\end{d_corollary}
\begin{proof}
Let $Y^*_{a,\tilde r}\ni \tilde a$ be a maximal self-stable family
in $\tilde Y$. Then, by Proposition \ref{6:p5.2} (i),
$[Y^*_{a,\tilde r}]_Z$ is maximal self-stable (in $\tilde Z$).
Since $\tilde a\in[Y^*_{a,\tilde r}]$, we have
$[Y^*_{a,r}]_Z=\tilde Z_{a,r}$. Hence, by \eqref{6:eq5.2},
$$
Y^*_{a,r}=[[Y^*_{a,r}]_Z]_Y=[\tilde Z_{a,\tilde r}]_Y=\tilde
Y_{a,\tilde r}.
$$
\end{proof}
Let $Y$ be a subspace of a metric space $(X,d)$. For $a\in Y$ and
$t>0$ we denote by
$$
S_t^Y=S^Y(a,t):=\{y\in Y:d(a,y)=t\}
$$
the sphere (in the subspace $Y$) with the center $a$ and the
radius $t$. Similarly for $a\in Z\subseteq X$ and $t>0$ define
$$
S_t^Z=S^Z(a,t):=\{z\in Z:d(a,z)=t\}.
$$
Write
\begin{equation}\label{6:eq5.6}
\varepsilon_a(t,Z,Y):=\sup_{z\in S_t^Z}\inf_{y\in Y}d(z,y)
\end{equation}
and
\begin{equation}\label{6:eq5.7}
\varepsilon_a(t)=\varepsilon_a(t,Z,Y)\wedge\varepsilon_a(t,Y,Z).
\end{equation}
\begin{d_theorem}\label{6:t5.4}
Let $Y$ and $Z$ be subspaces of a metric space $(X,d)$ and let
$a\in Y\cap Z$. Then $Y$ and $Z$ are strongly tangent equivalent
at the point $a$ if and only if the equality
\begin{equation}\label{6:eq5.8}
\lim_{t\to0}\frac{\varepsilon_a(t)}t=0
\end{equation}
holds.
\end{d_theorem}
\begin{proof}
Suppose that limit relation \eqref{6:eq5.8} holds. Let $\tilde
r=\{r_n\}_{n\in\mathbb N}$ be a normalizing  sequence and $\tilde
z=\{z_n\}_{\in\mathbb N}\in\tilde Z$ be a sequence with a finite
limit
$$
\tilde d_{\tilde r}(\tilde a,\tilde z)=\lim_{n\to
\infty}\frac{d(a,z_n)}{r_n}.
$$
To find $\tilde y=\{y_n\}\in\tilde Y$ such that
\begin{equation}\label{6:eq5.9}
\tilde d_{\tilde r}(\tilde y,\tilde z)=0
\end{equation}
note that we can take $\tilde y=\tilde a$ if $\tilde d_{\tilde
r}(\tilde a,\tilde z)=0$. Hence, without loss of generality, we
suppose
\begin{equation}\label{6:eq5.10}
0<\tilde d_{\tilde r}(\tilde a,\tilde z)=\lim_{n\to
\infty}\frac{d(z_n,a)}{r_n}<\infty.
\end{equation}
It follows from \eqref{6:eq5.7} and \eqref{6:eq5.8} that
\begin{equation}\label{6:eq5.11}
\lim_{t\to0}\frac1t(\varepsilon_a(t,Z,Y))=0.
\end{equation}
Inequalities \eqref{6:eq5.10} imply that there is $n_0\in\mathbb
N$ such that $d(z_n,a)>0$ if $n\geq n_0$. Write for every
$n\in\mathbb N$
\begin{equation}\label{6:eq5.12}
t_n:=\begin{cases} 1&\text{if }n<n_0\\
d(z_n,a)&\text{if }n\geq n_0.
\end{cases}
\end{equation}
The definition of $\varepsilon_a(t,Z,Y)$ implies that for every
$n\in\mathbb N$ there is $y_n\in Y$ with
\begin{equation}\label{6:eq5.13}
d(z_n,y_n)\leq\varepsilon_a(t_n,Z,Y)+t_n^2.
\end{equation}
 Put $\tilde
y=\{y_n\}_{n\in\mathbb N}$ where $y_n$ are points in $Y$ for which
\eqref{6:eq5.13} holds. Now using
\eqref{6:eq5.10}--\eqref{6:eq5.12} we obtain
\begin{multline*}
\limsup_{n\to\infty}\frac{d(z_n,y_n)}{r_n}\leq\lim_{n\to\infty}\frac{d(z_n,
a)}{r_n}\limsup_{n\to\infty}\frac{d(z_n,y_n)}{t_n}\\\leq \tilde
d_{\tilde r}(\tilde z,\tilde
a)\limsup_{n\to\infty}\frac{\varepsilon_a(t_n,Z,Y)+t_n^2}{t_n}=0.
\end{multline*}
Consequently $\lim_{n\to\infty}\frac{d(z_n,y_n)}{r_n}=0$, i.e.,
\eqref{6:eq5.9} holds. Similarly we can prove that for every
$\tilde y\in\tilde Y$ with a finite $\tilde d_{\tilde r}(\tilde
y,\tilde a)$ there is $\tilde z\in\tilde Z$ such that $\tilde
d_{\tilde r}(\tilde z,\tilde y)=0$. Hence if \eqref{6:eq5.8}
holds, then $Y$ and $Z$ are strongly tangent equivalent at the
point $a$.

Suppose now that \eqref{6:eq5.8} does not hold. More precisely, we
shall assume that
$$
\limsup_{n\to\infty}\frac{\varepsilon_a(t,Z,Y)}{t}>0.
$$
Then there is a sequence $\tilde t$ of positive numbers $t_n$ with
$\lim_{n\to\infty}t_n=0$ and there is $c>0$ such that for every
$n\in\mathbb N$ there exists $z_n\in S^Z(a,t_n)$ for which
\begin{equation}\label{6:eq5.14}
\inf_{y\in Y}d(z_n,y)\geq ct_n=cd(a,z_n).
\end{equation}
Let us denote by $\tilde z$ the sequence of points $z_n\in Z$
which satisfy \eqref{6:eq5.14}. Take the sequence $\tilde
t=\{t_n\}_{n\in\mathbb N}$ as a normalizing sequence. Then, by
\eqref{6:eq5.14}, we obtain
$$
\limsup_{n\to\infty}\frac{d(z_n,y_n)}{t_n}\geq c>0
$$
for every $\tilde y=\{y_n\}_{n\in\mathbb N}\in\tilde Y$.
Consequently $Y$ and $Z$ are not strongly tangent equivalent at
the point $a$.
\end{proof}
Consider now the case where $Z=X$. Let $(X, d)$ be a metric space
and let $a\in Y\subseteq X$. If $\tilde X_{a, \tilde r}$ is a
maximal self-stable family of sequences $\tilde x\in \tilde X$ and
if $\tilde Y_{a, \tilde r}=\tilde Y\cap\tilde X_{a,\tilde r}$,
then it is obvious that $\tilde Y_{a,\tilde r}$ is also maximal
self-stable (in $\tilde Y$) and that there is a unique isometric
embedding $E_{m_Y}: \Omega_{a, \tilde r}^Y \to \Omega_{a, \tilde
r}$ such that the following diagram
\begin{equation}\label{6:eq5.15}
\begin{array}{ccc}
\tilde Y_{a, \tilde r}&\xrightarrow{\ \ \ \ \mbox{in}\ \ \ \ }&
\tilde X_{a, \tilde r}
\\
\negthinspace \negthinspace\negthinspace\negthinspace\! \!\! \!\!
\! p_Y\Bigg\downarrow & &\negthinspace\! \!\!\Bigg\downarrow p
\\
\Omega_{a, \tilde r}^Y&\xrightarrow {\ \ \ E_{m_Y}\ \ \ }&
\Omega_{a, \tilde r}
\end{array}
\end{equation}
is commutative. Here $\Omega_{a, \tilde r}$ is a pretangent space
correspondent to $\tilde X_{a,  r}$, $\Omega_{a, \tilde r}^Y$ is a
metric identification of $\tilde Y_{a, \tilde r}$, $p_Y$ and $p$
are appropriate metric identification maps and ${\rm in}(\tilde
y)=\tilde y$ for all $\tilde y \in \tilde Y_{a, \tilde r}$.
\begin{d_corollary}\label{6:c5.5}
Let $(X, d)$ be a metric space, let $Y$ be a subspace of $X$ and
let $a\in Y$. The following conditions are equivalent.
\begin{itemize}
\item[$(i)$] An embedding $E_{m_Y} : \Omega_{a, \tilde
    r}^Y \to \Omega_{a, \tilde r}$ is an isometry for every
normalizing sequence $\tilde r$ and for all maximal self-stable
families $\tilde Y_{a, \tilde r}$ and $\tilde X_{a,
    \tilde r}$ with $\tilde Y_{a, \tilde r}\subseteq \tilde
    X_{a, \tilde r}$.

\item[$(ii)$] The equality
    $$
    \lim_{t\to 0} \frac{\varepsilon_a(t,X,Y)}{t}=0
    $$
    holds.

\item[$(iii)$] $X$ and $Y$ are strongly tangent equivalent at the
point $a$.
\end{itemize}
\end{d_corollary}
\begin{proof}
The equivalence $(ii)\Leftrightarrow(iii)$ follows from Theorem
\ref{6:t5.4}. To prove $(iii)\Rightarrow(i)$ note that
$$
[\tilde X_{a,\tilde r}]_Y=\tilde Y\cap\tilde X_{a,\tilde r}
$$
for every maximal self-stable $\tilde X_{a,\tilde r}$.
Consequently $(iii)$ implies $(i)$ because mapping \eqref{6:eq5.3}
is an isometry.

Now suppose that the mappings $E_{m_Y}$ from \eqref{6:eq5.15} are
isometries for all $\tilde X_{a,\tilde r}$. To prove
$(i)\Rightarrow(iii)$ it is sufficient to show that for every
maximal self-stable $\tilde X_{a,\tilde r}$ and every $\tilde
x_0\in\tilde X_{a,\tilde r}$ there is $\tilde y_0\in\tilde Y$ such
that
\begin{equation}\label{6:eq6.16}
\tilde d_{\tilde r}(\tilde x_0,\tilde y_0)=0.
\end{equation}
Let $\alpha=p(\tilde x_0)$. Since $E_{m_Y}$ is an isometry,
$E_{m_Y}$ is surjective. Thus $E^{-1}_{m_Y}(\alpha)\ne\emptyset$.
The last condition and
$$
p_Y^{-1}(E^{-1}_{m_Y}(\alpha))\ne\emptyset
$$
are equivalent because $p_Y$ also is surjective. Since
$$E_{m_Y}\circ p_Y=p\circ\text{in}$$ and $p^{-1}(\alpha)=\{\tilde
x\in\tilde X_{a,\tilde r}:\tilde d_{\tilde r}(\tilde x,\tilde
x_0)=0\},$ we have
\begin{multline*}
\emptyset\ne
p^{-1}_Y(E^{-1}_{m_Y}(\alpha))=\text{in}^{-1}(p^{-1}(\alpha)) \\=
\text{in}^{-1}(\{\tilde x\in\tilde X_{a,\tilde r}:\tilde d_{\tilde
r}(\tilde x_0,\tilde x)=0\})\\=\tilde Y\cap\{\tilde x\in\tilde
X:\tilde X:\tilde d_{\tilde r}(\tilde x_0,\tilde x)=0\},
\end{multline*}
that implies \eqref{6:eq6.16} with some $\tilde y_0\in\tilde Y$.
\end{proof}
Obviously, condition $(ii)$ of Corollary \ref{6:c5.5} holds if $Y$
is a dense subset of $X$. Therefore we have the following.

\begin{d_corollary}\label{6:c5.6}
Let $(X, d)$ be a metric space and let $Y$ be a dense subspace of
$X$. Then $X$ and $Y$ are strongly tangent equivalent at all
points $a\in Y$, in particular, the pretangent spaces to $X$ and
to $Y$ are pairwise isometric for all normalizing sequences at
every point $a\in Y$.
\end{d_corollary}

Consider now some examples.

The following result was proved in \cite{DAK1}. Let $X=\mathbb R$
or $X=\mathbb C$ or $X=\mathbb R^+=[0,\infty)$ and let
$$
d(x,y)=|x-y|
$$
for all $x,y\in X$.
\begin{d_proposition}\label{p:3.2}
Each pretangent space $\Omega_{0,\tilde r}$ (to $X$ at the point
$0$) is tangent and isometric to $(X,d)$ for every normalizing
sequence $\tilde r$.
\end{d_proposition}
Using Theorem \ref{6:t5.4} and Proposition \ref{p:3.2} we can
easily obtain future examples of tangent spaces to some subspaces
of the Euclidean space $E^n$. The first example will be examined
in details.
\begin{d_example}\label{e3.1}
Let $F:[0,1]\to E^n,\ n\geq 2$, be a simple closed curve in the
Euclidean space $E^n$, i.e., $F$ is continuous and $F(0)=F(1)$ and
$$
F(t_1)\ne F(t_2)
$$
for every two distinct points $t_1,t_2\in[0,1]$ with
$|t_2-t_1|\ne1$. We can write $F$ in the  coordinate form
$$
F(t)=(f_1(t),\dots,f_n(t)),\qquad t\in[0,1].
$$
Suppose that all functions $f_i,\ 1\leq i\leq n$, are
differentiable at a point $t_0\in(0,1)$ and
$$
|F'(t_0)|=\bigg|\sum_{i=1}^n(f'_i(t_0))^2\bigg|^{\frac12}\ne0.
$$
(In the case $t_0=0$ or $t_0=1$ we must use the one-sided
derivatives.) We claim that each pretangent space to the subspace
$Y=F([0,1])\subseteq E^n$ at the point $a=F(t_0)$ is tangent and
isometric to $\mathbb R$ for every normalizing sequence $\tilde
r$. Indeed, by Proposition \ref{6:p5.2} and \ref{p:3.2}, it is
sufficient to show that $Y$ is strongly tangent equivalent to the
straight line
$$
Z=\{(z_1(t),\dots,z_n(t)):(z_1(t),\dots,z_n(t))=F'(t_0)(t-t_0)+F(t_0),\
t\in\mathbb R\}
$$
at the point $a=F(t_0)$.

The classical definition of the differentiability of real
functions shows that limit relation \eqref{6:eq5.5} holds with
these $Y$ and $Z$. Hence, by Theorem~\ref{6:t5.4}, $Y$ and $Z$ are
strongly tangent equivalent at the point $a=F(t_0)$.
\end{d_example}

\begin{d_example}\label{e3.2}
Let $f_i:[-1,1]\to\mathbb R,\ i=1,\dots,n$, be functions such that
$f_1(0)=\dots=f_n(0)=c$ where $c\in\mathbb R$ is a constant.
Suppose all $f_i$ have a common finite derivative $b$ at the point
$0$, $f_1'(0)=\dots=f_n'(0)=b$. Write
$$
a=(0,c)\quad\text{and}\quad
X=\bigcup_{i=1}^n\{(t,f_i(t)):t\in[-1,1]\},
$$
i.e., $X$ is an union of the graphs of the functions $f_i$. Let us
consider $X$ as a subspace of the Euclidean plane $E^2$. Then each
pretangent space $\tilde X_{a,\tilde r}$ to the space $X$ at the
point $a$ is tangent and isometric to $\mathbb R$.
\end{d_example}

\begin{d_example}\label{e3.3}
Let $f_1,f_2$ be two functions from the precedent example. Put
$$
X=\{(x,y):f_1(x)\wedge f_2(x)\leq y\leq f_1(x)\vee f_2(x),\
x\in[-1,1]\},
$$
i.e., $X$ is the set of points of the plane which lie between the
graphs of the functions $f_1$ and $f_2$. Then each pretangent
space $\tilde X_{a,\tilde r}$ to $X$ at $a=(0,0)$ is tangent and
isometric to $\mathbb R$.
\end{d_example}
\begin{d_example}\label{e3.4}
Let $\alpha$ be a positive real number. Write
$$
X=\{(x,y,z)\in E^3:\sqrt{y^2+z^2}\leq x^{1+\alpha},\ x\in\mathbb
R^+\},
$$
i.e., $X$ can be obtained by the rotation of the plane figure
$\{(x,y)\in E^2:0\leq y\leq x^{1+\alpha},\ x\in\mathbb R^+\}$
around the real axis.  Then each pretangent space $\tilde
X_{a,\tilde r}$ to $X$ at the point $a=(0,0,0)$ is tangent and
isometric to $\mathbb R^+$.
\end{d_example}

\begin{d_example}\label{e3.5}
Let $U\subseteq \mathbb C$ be an open set, let $F:U\to E^n,\
n\geq2$, be an one-to-one continuous function,
$$
F(x,y)=(f_1(x,y),\dots,f_n(x,y)),\quad (x,y)\in U,
$$
and let $(x_0,y_0)$ be a marked point of $U$. Suppose that all
$f_i$ are  differenti\-able at the point $(x_0,y_0)$ and that the
rank of the Jakobian matrix of $F$ equals two at this point. Write
$$
X=F(U),\quad a=F(x_0,y_0).
$$
Consider the parametrized surface $X$ as a subspace of $E^n$. Then
every pretangent space $\Omega_{a,\tilde r}^X$ is tangent and
isometric to $\mathbb C$.
\end{d_example}

{\bf Acknowledgments.} The author  thanks the Department of
Mathe\-matics and Statistics of the University of Helsinki for the
comfortable setting in the May--June 2008 when he began to work
with the initial version of this paper. This work  also was
partially supported by the Ukrainian State Foundation for Basic
Researches, Grant $\Phi$ 25.1/055.

\end{document}